\documentclass[a4paper,10pt]{article}
\usepackage[utf8x]{inputenc}
\usepackage[russian]{babel}
\usepackage{amssymb}

\title{Симплектическая редукция и  лагранжевы подмногообразия  в ${\rm Gr}(1, n)$}
\author{ Николай Тюрин\\
ЛТФ ОИЯИ (Дубна) и МИ РАН (Москва)}

\date{}

\begin{document}

\maketitle

\begin{abstract} В работе построены новые примеры лагранжевых подмногообразий комплексного грассманиана ${\rm Gr}(1, n)$, снабженного стандартной кэлеровой формой. Схема построения
исходит из двух фактов: во - первых, мы предлагаем естественное соответствие между лагранжевыми подмногообразиями в симплектическом многообразии, являющимся результатом симплектической редукции, и лагранжевыми подмногообразиями большого симплектического многообразия с гамильтоновым действием группы, к которому применяется эта редукция; во- вторых, мы показываем,
что при некотором подборе порождающих действия ${\rm T}^k$ при $k = 2, ..., n-1$, на ${\rm Gr}(1, n)$ и подходящих значениях отображений моментов имеется изоморфизм ${\rm Gr}(1, n) \slash  \slash {\rm T}^k \cong {\rm tot}(\mathbb{P}(\tau) \times ... \mathbb{P}(\tau) \to {\rm Gr}(1, n-k))$, где справа стоит тотальное пространство прямого произведения $k$ копий проективизации тавтологического расслоения $\tau \to {\rm Gr}(1, n-k)$. Комбинируя эти два факта мы получаем нижнюю оценку
на число топологически различных гладких лагранжевых подмногообразий в исходном грассманиане ${\rm Gr}(1, n)$.
\end{abstract}

\section*{Введение}

В статье  [1]  А. Миронов предложил новую конструкцию  подмногообразий в  $\mathbb{C}^n$ и $\mathbb{C} \mathbb{P}^n$, удовлетворяющих условию лагранжевости относительно
  кэлеровых форм, индуцируемых  постоянной метрикой и метрикой Фубини - Штуди соответственно. В своей работе он исходил из задачи построения минимальных и гамильтоново минимальных лагранжевых подмногообразий и лагранжевых погружений (то есть лагранжевых подмногообразий, допускающих самопересечения). В  работе [2] было показано,    что конструкция Миронова может быть обобщена на случай произвольного кэлерова многообразия с действием тора, однако в отличие от исходной работы там не обсуждалась тема
   минимальности получаемых лагранжевых циклов. В настоящей работе мы используем простое обобщение этих конструкции  на еще более общий случай, когда на симплектическом многообразии имеется гамильтоново действие произвольной  группы Ли.

   Пусть $(M, \omega)$ --- компактное симплектическое многообразие, на котором гамильтоново действует группа Ли $G$ с отображением моментов $\mu: M \to {\frak g}^*$
   в двойственную алгебру Ли. Пусть $\mu^{-1} (t) \subset M$ множество уровня, соответствующее центральному элементу $t  \in {\frak g}^*$, на котором $G$ действует свободно,
   так что имеется фактор - многообразие $N_t = \mu^{-1}(t)/G$ с соответствующей симплектической формой $\omega_t$, подлежащее главному расслоению $\pi_t: \mu^{-1}(t) \to N_t$. Тогда имеется следующая

   {\bf Теорема 1.} {\it  Пусть $S_0 \subset N_t$ --- лагранжево вложение (погружение). Тогда его прообраз $S = \pi^{-1}_t (S_0) \subset \mu^{-1}(t) \subset M$ есть лагранжево вложение
   (погружение) в $(M, \omega)$}.

   В самом деле, по определению симплектической редукции  $M // G$, см. [4], касательное пространство $T_p S$ составлено из двух трансверсальных компонент: ${\rm ker} ({\rm d} \pi_t|_p) \oplus    T_{\pi_t (p)}S_0$, первая из которых естественно изоморфна ${\frak g}$ и содержится в ${\rm ker} \omega_p|_{T_p \mu^{-1}(t)}$. Отсюда из определения формы $\omega_t$ и лагранжевости
   $S_0$ следует утверждение теоремы.

   {\bf Замечание 0.} Нетрудно видеть, что если на исходном многообразии $M$ имеется произвольное лагранжево подмногообразие $K \subset M$, трансверсальное действию группы $G$, то
    этим определено лагранжево погружение $N_t \supset S_0 = \pi_t(\mu^{-1}(t) \cap K)$, применяя к которому Теорему 1 мы получим новое лагранжево подмногообразие
    (переход от $K$ к этому последнему мы называем перестройкой лагранжевых циклов). Отсюда видно, что  основные утверждения работы [2] следуют из Теоремы 1, если в качестве $G$ рассмотреть тор ${\rm T}^k$, а в качестве $K$ --- вещественную часть $M$ относительно действия антиголоморфной инволюции.
   Поэтому естественно понимать Теорему 1 как обобщение конструкций из работ [1], [2].

В работе [2] было показано, что комплексный грассманиан ${\rm Gr}(1, n)$ прямых в проективном пространстве $\mathbb{C} \mathbb{P}^n$, снабженный стандартной кэлеровой формой, порождаемой плюккеровым вложением, допускает гамильтоново действие тора ${\rm T}^n$, сохраняющего и комплексную структуру. В этой же работе была высказана гипотеза о том, что  ${\rm Gr}(1, n)$
допускает  не менее $n$ различных топологических типов, реализуемых гладкими  лагранжевыми вложениями. Данная гипотеза  исходила  из того соображения, что для каждого $k$ от 0 до $n$ мы можем применить
обобщенную конструкцию из [2] и получить {\it a priori} топологически различные лагранжевы подмногообразия.  Теорема 1 подсказывает описание таких обобщенных циклов Миронова различной степени однородности в  ${\rm Gr}(1, n)$, так что оценка, приведенная в гипотезе, оказывается превзойденной.
А именно, ниже в настоящей работе мы доказываем, что

{\bf Теорема 2.} {\it Для любого $k \in (2, ..., n-1)$  на грассманиане ${\rm Gr}(1, n)$ существуют  функции $\tilde \mu_0, ..., \tilde \mu_{k-1}$, индуцирующие гамильтоново действие ${\rm T}^k$, и набор значений  $c_0, ..., c_{k-1} \in \mathbb{R}_{+}$, таких что
результатом соответствующей симплектической редукции является ${\rm Gr}(1, n) \slash \slash {\rm T}^k \cong M_{n-k} = {\rm tot}(\mathbb{P}(\tau) \times ... \mathbb{P}(\tau) \to {\rm Gr}(1, n-k))$, где справа стоит тотальное пространство прямого произведения $k$ копий проективизации тавтологического расслоения $\tau \to {\rm Gr}(1, n-k)$.}

Более того, в результате редукции остальные образующие действия ${\rm T}^n$ на ${\rm Gr}(1, n)$ не пропадают, а порождают естественное гамильтоново действие тора  ${\rm T}^{n-k}$ на  многообразии $M_{n-k}$. Явно описать это действие возможно следующим образом. Реализуем $M_{n-k}$ как подмногообразие в $k+1$ - ом прямом произведении
$\mathbb{C} \mathbb{P}^{n-k} \times ... \times \mathbb{C} \mathbb{P}^{n-k} \times {\rm Gr}(1, n-k)$, задаваемое пересечением $k$ "частичных циклов инциденции" ${\cal U}_i $ в $k$- ом прямом произведении $ \mathbb{C} \mathbb{P}^{n-k} \times ... \times \mathbb{C} \mathbb{P}^{n-k} \times {\rm Gr}(1, n-k)$,
составленных из пар (точка, прямая), где для ${\cal U}_i$ точка берется из $i$ - го прямого сомножителя. При этом на каждом из $k+1$ прямом слагаемом
после отождествления первых $k$ диагонально действует тор ${\rm T}^{n-k}$, очевидным образом сохраняющий подмногообразие $M_{n-k}$. Поэтому симплектическая редукция ${\rm Gr}(1, n) \slash \slash {\rm T}^k$ может быть представлена как симплектическая редукция $M_{n-i} \slash \slash {\rm T}^{k-i}$ для $i =0, ..., k$.

С другой стороны, для любого $k$ действие тора $T^k$, описанное выше, сохраняет всю кэлерову структуру, откуда  многообразие $M_{n-k}$ естественным образом наделено кэлеровой структурой,
откуда следует что вещественная часть $M_{n-k}^{\mathbb{R}} \subset M_{n-k}$ является изотропным относительно кэлеровой формы. Поскольку $M_{n-k}^{\mathbb{R}} \cong{\rm tot}(\mathbb{P}(\tau_{\mathbb{R}}) \times ... \mathbb{P}(\tau_{\mathbb{R}}) \to {\rm Gr}_{\mathbb{R}} (1, n-k))$, где справа стоит тотальное пространство прямого произведения $k$ копий проективизации тавтологического расслоения $\tau_{\mathbb{R}} \to {\rm Gr}_{\mathbb{R}}(1, n-k)$, то его вещественная  размерность в точности равна комплексной размерности
 всего $M_{n-k}$, следовательно оно является лагранжевым в $M_{n-k}$. Заметим, что топологически $M^{\mathbb{R}}_{n-k}$ есть расслоение на $k$ - мерный тор над ${\rm Gr}_{\mathbb{R}}(1, n-k)$, поэтому из Теоремы 1 мы получаем  семейство лагранжевых подмногообразий $\{ S_k, k =2, ..., n-1 \}$, таких что каждое $S_k$ изоморфно расслоению со слоем $T^{2k}$ над ${\rm Gr}_{\mathbb{R}}(1, n-k)$. Нетрудно видеть, что для разных $k$ соответствующие $S_k$ топологически не эквивалентны,  что очевидным образом влечет положительный ответ для Гипотезы, сформулированной в конце работы [2]. Как будет показано ниже, дополнительный анализ наших конструкций доставляет еще не менее чем $[\frac{n}{2}] \cdot [\frac{n-1}{2}]$
  новых топологических типов,  что существенно  увеличивает оценку на различные топологические типы гладких лагранжевых подмногообразий в ${\rm Gr} (1, n)$.

Поскольку обсуждение условий и доказательство Теоремы 1 уместилось выше, в следующих параграфах мы сразу переходим к доказательству Теоремы 2: сначала разберем
геометрические и аналитические аспекты для случая ${\rm Gr}(1, 3)$, а затем, базируясь на полученном материале, разберем общий случай ${\rm Gr}(1, n)$.

\section{Грассманиан ${\rm Gr}(1, 3)$}

Рассмотрим проективное пространство $\mathbb{C} \mathbb{P}^3$ с фиксированными координатами $[z_0: ... :z_3]$ и согласованной стандартной кэлеровой формой $\omega_{FS}$
соответствующей метрики Фубини - Штуди. Имеем тогда набор отображений моментов
$$
\mu_i = \frac{\vert z_i \vert^2}{\sum_{i=0}^3 \vert z_i \vert^2}, \quad i = 0, ..., 3,
$$
из которых любые три индуцируют гамильтоново действие тора ${\rm T}^3$. Каждое из отображений $\mu_i$ индуцирует действие на множестве всех проективных прямых, составляющих
многообразие ${\rm Gr}(1, 3)$, откуда получаем индуцированное действие ${\rm T}^3$ на последнем многообразии, которое является гамильтоновым относительно стандартной кэлеровой формы (все явные формулы будут приведены ниже, в аналитической части конструкции). Рассмотрим два отображения моментов $\mu_0, \mu_1$ и исследуем соответствующие им отображения моментов $\tilde \mu_i$ на грассманиане
${\rm Gr}(1, 3)$. Обозначим как $p_i, i =0, ..., 3,$ точки в $\mathbb{C} \mathbb{P}^3$, соответствующие базисным ортонормированным векторам $e_i$ в $\mathbb{C}^4$, проективизацией которого
является наше проективное пространство. Тогда (см [3]) функции $\tilde \mu_i$ представляются как
$$
\tilde \mu_i (l) = \max_{v \in V} \frac{\vert <e_i; v> \vert^2}{\vert v \vert^2},
\eqno (1)
$$
где $V \subset \mathbb{C}^4$ --- двумерное продпространство, соответствующее проективной прямой $l \subset \mathbb{C} \mathbb{P}^3$. Отсюда видно, что функция
$\tilde \mu_i$ принимает значение на отрезке $[0; 1]$, ее критическими значениями являются $0$ и $1$, которые достигаются в случае если $l \subset \langle p_{j_1}, p_{j_2}, p_{j_3} \rangle, j_k \neq i,$
или $p_i \in l$ соответственно (под $\langle  ... \rangle$ в этой работе  мы будем обозначать проективную оболочку объектов, содержащихся в этих скобках). Таким образом, критическим множеством для функции $\tilde \mu_0$
 будет подмножество прямых, проходящих через точку $p_0$, в объединении с множеством прямых, содержащихся в плоскости $\langle p_1, p_2, p_3 \rangle$, то есть соответствующие $\alpha$- и $\beta$-
 плоскости в ${\rm Gr}(1, 3)$. Кроме того, нетрудно видеть, что если прямая $l \subset \mathbb{C} \mathbb{P}^3$ имеет нетривиальное пересечение с прямой $ \langle p_0, p_1 \rangle$, то $\tilde \mu_0(l) + \tilde \mu_1(l) \geq 1$. В самом деле, если пересечение нетривиально, то в соответствующем $V$ найдется вектор $v \in V$, представляющийся как $v = x_0 e_0 + x_1 e_1$, откуда
 из формулы (1) получаем, что сумма максимумов должна быть не меньше 1.

 Рассмотрим открытое подмножество ${\rm Gr}^0(1, 3) \subset {\rm Gr}(1, 3)$, определяемое условиями
 $$
 {\rm Gr}^0(1, 3) = \{l \in {\rm Gr}(1, 3) \quad \vert \quad \tilde \mu_0(l) > 0, \tilde \mu_1(l) >0, \tilde \mu_0(l) + \tilde \mu_1(l) <1 \}.
 $$

 Тогда имеем

 {\bf Предложение 1.} {\it Открытое многообразие ${\rm Gr}^0(1, 3)$ расслоено над прямым произведением $\langle p_2, p_3 \rangle \times \langle p_2, p_3 \rangle$ так, что
 гамильтоново действие ${\rm T}^2$, натянутое на $\tilde \mu_0$ и $\tilde \mu_1$, сохраняет базу и слои.}

{\bf Доказательство.} Рассмотрим  произвольную прямую $l \subset \mathbb{C} \mathbb{P}^3$, которая  не пересекает $\langle p_0, p_1 \rangle$ и не лежит в подпространствах $\langle p_i, p_2, p_3 \rangle, i =0, 1$.
 Тогда корректно определены точки $s_0, s_1 \in \langle p_2, p_3 \rangle$, получаемые как пересечения плоскостей, натянутых на $p_1$ и $l$ или $p_0$ и $l$ соответственно,  и прямой $\langle p_2, p_3 \rangle$. В самом деле,
 плоскость $\langle p_0, l \rangle \subset \mathbb{C} \mathbb{P}^3$ не может содержать всю прямую $\langle p_2, p_3 \rangle$, поскольку тогда $l$ обязана была бы лежать в плоскости $\langle p_0, p_2, p_3 \rangle$,  что исключено по условию; те же аргументы применимы и к плоскости $\langle p_1, l \rangle$. Таким образом, прямая $l$ соответствует паре точек $s_0, s_1 \in \langle p_2, p_3 \rangle$ (возможно, совпадающих).
 При этом соотвествие не симметрично (но симметрия очевидным образом реализуется перестановкой двух отображений моментов, что в свою очередь соответствует автоморфизму исходного
 проективного пространства $\mathbb{C} \mathbb{P}^3$, реализуемому перестановкой однородных координат).

 С другой стороны, как восстанавливается прямая $l$ по паре точек $(s_0, s_1)$? Рассмотрим два пучка плоскостей, проходящих через прямые $\langle p_0, s_0 \rangle$ и $\langle p_1, s_1 \rangle$ соответственно.
  Из каждого пучка  удалим пару плоскостей, содержащих прямые $\langle p_0, p_1 \rangle, \langle p_2, p_3 \rangle$. Тогда нетрудно видеть, что ${\rm Gr}^0(1, 3)$ естественным образом вкладывается в расслоение
 над прямым произведением $\langle p_2, p_3 \rangle \times \langle p_2, p_3 \rangle$ и со слоем $\mathbb{C}^* \times \mathbb{C}^*$. В самом деле, по построению в пучках плоскостей обязаны присутствовать  плоскости $\langle p_0, l \rangle$ и $ \langle p_1, l \rangle $, каждая из которых лежит в соответствующем пучке, а их пересечение в точности дает прямую $l$. При этом грубая оценка показывает,
 что образ вложения ${\rm Gr}^0(1, 3)$ в расслоение равномерно распределен: над каждой точкой $\langle p_2, p_3 \rangle  \times \langle p_2, p_3 \rangle$ имеем открытое четырехмерное подмножество в слое.
 Это означает, что при выборе конкретных положительных значений $c_0, c_1 \in \mathbb{R}_+, c_0 + c_1 < 1,$ совместное множество уровня $N_{c_0, c_1}
 = \{ \tilde \mu_0 = c_0, \tilde \mu_1 = c_1 \} \subset {\rm Gr}(1, 3)$
 расслоено над прямым произведением $\langle p_2, p_3 \rangle \times \langle p_2, p_3 \rangle$ со слоем ${\rm T}^2$ (формульно этот факт подтверждается ниже).

 {\bf Замечание 1.} При  этом нетрудно видеть, что  расслоение со слоем $\mathbb{C}^* \times \mathbb{C}^*$ над прямым произведением $\langle p_2, p_3 \rangle \times \langle p_2, p_3 \rangle$  не является топологически тривиальным:  оно естественно вложено в расслоение со слоем $\mathbb{C} \mathbb{P}^1 \times \mathbb{C} \mathbb{P}^1$, причем каждое из этих  проективных расслоений над проективной прямой $\langle p_2, p_3 \rangle$  допускает следующую  реализацию, позволяющую определить их топологический тип. Поскольку ответы в обоих случаях совпадают, рассмотрим пучок плоскостей, проходящих через $\langle p_0, s_0 \rangle$.  Рассмотрим пересечение
 каждой плоскости с проективной оболочкой $\langle p_1, p_2, p_3 \rangle$, что сопоставляет соответствующей плоскости прямую, проходящую через $s_0$. Тогда как расслоение над $\langle p_2, p_3 \rangle$ получаем в качестве слоя пучок прямых, проходящих через точку $s_0$, что при глобализации дает тотальное пространство, изоморфное проективизации расслоения ${\cal O} \oplus {\cal O}(-1)$  над проективной прямой $\langle p_2, p_3 \rangle$. Это расслоение допускает пару глобальных сечений, соответствующих плоскостям $\langle p_0, s_0,  p_1 \rangle$ и $\langle p_0, p_2, p_3 \rangle$, при этом выбор значений $\tilde \mu_0 = c_0, \tilde \mu_1 = c_1$ для $c_i >0, c_0 +c_1 < 1$ выделяет в слоях пару окружностей, так что тотальное пространство каждого из расслоений
 изоморфно $S^3$ (тотальному пространству расслоения Хопфа), откуда
 $$
{\rm Gr}^0(1, 3) \supset \{ \tilde \mu_0^{-1}(c_0) \cap \tilde \mu_1^{-1}(c_1) \} \cong S^3 \times S^3 \to \langle p_2, p_3 \rangle \times  \langle p_2, p_3 \rangle.
 $$
 При этом подчеркнем, что для расслоения Хопфа $\pi: S^3 \to S^2$ ограничение на произвольную гладкую петлю $\gamma \subset S^2$ дает гладкий двумерный тор $\pi^{-1}(\gamma) \subset S^3$.

 Обратим внимание теперь на следующий факт: гамильтоново $U(1)$ - действие, порождаемое $\tilde \mu_0$ на ${\rm Gr}(1, 3)$, оставляет на месте базу расслоения и действует на слоях,
 сохраняя структуру прямого произведения. В самом деле, пусть прямая $l$ есть пересечение плоскостей $\pi_0$ и $\pi_1$ из пучков плоскостей, содержащих прямые
 $\langle p_0, s_0 \rangle$ и $ \langle p_1, s_1 \rangle$ соответственно. Гамильтоново действие каждого из $\tilde \mu_i$ оставляет неподвижными все точки $p_i, s_j$, поскольку эти точки содержатся в критических
 множествах $\tilde \mu_i$. Отсюда следует, что прямые $\langle p_i, s_i \rangle$ переходят в себя при этом действии, и следовательно пучки плоскостей также переходят в себя. Отсюда видно,
 что торическое действие сохраняет  базу $\langle p_2, p_3 \rangle  \times  \langle p_2, p_3 \rangle$ и структуру расслоения, действуя послойно. Этим завершается доказательство Предложения 1.

 Заметим, что для вывода из Предложения 1 утверждения Теоремы 2 для случая ${\rm Gr}(1, 3)$ нам осталось показать, что совместное множество уровня функций $\tilde \mu_i$ есть
 расслоение со слоем ${\rm T}^2$. Проведем аналитические вычисления и выведем явные формулы.

 В плюккеровых переменных $w_{ij}$ на проективном пространстве $\mathbb{P}(\Lambda^2 \mathbb{C}^4)$ где $0 \leq i < j \leq 3$
грассманиан ${\rm Gr}(1, 3)$ реализуется квадрикой $Q = \{ w_{01}w_{23} - w_{02} w_{13} + w_{03} w_{12} = 0 \} \subset \mathbb{C} \mathbb{P}^5$; наши выделенные отображения моментов имеют вид
$$
 \tilde \mu_0 = \frac{\vert w_{01} \vert^2 + \vert w_{02} \vert^2 + \vert w_{03} \vert^2}{\sum_{i<j} \vert w_{ij} \vert^2},
 \quad \tilde \mu_1 =  \frac{\vert w_{01} \vert^2 + \vert w_{12} \vert^2 + \vert w_{13} \vert^2}{\sum_{i<j} \vert w_{ij} \vert^2},
 \eqno (2)
 $$
 (см. [2]). Легко видеть, что соответствующее гамильтоново действие сохраняет квадрику $Q$, но кроме того на квадрике $Q$ имеется два пучка дивизоров, каждый из которых
 инвариантен относительно торического действия. В самом деле, рассмотрим пучки
 $$
 Q_0 (\alpha) = \{ \alpha_0 w_{02} + \alpha_1 w_{03} = 0 \} \cap Q, \quad Q_1(\beta) =  \{\beta_0 w_{12} + \beta_1 w_{13} \} \cap Q.
  $$
  Так как каждая из $\tilde \mu_i$ тривиально действует на $Q_j$ если $i \neq j,$, и при этом действует одновременным умножением
  на $e^{it}$ каждой координаты, входящей в определение $Q_i$, то соответствующие пучки инвариантны относительно всего торического действия ${\rm T}^2$. По определению плюккерова вложения
  пересечение $N(\alpha, \beta) = Q_0(\alpha) \cap Q_1 (\beta) \subset {\rm Gr}(1, 3)$ состоит из прямых, имеющих нетривиальное пересечение с двумя соответствующими  прямыми из пучков,
  натянутых в $\mathbb{C} \mathbb{P}^3$ на пары $( \langle p_0, p_2 \rangle, \langle p_0, p_3 \rangle)$ и $(\langle p_1, p_2 \rangle, \langle p_1, p_3 \rangle)$. Каждая прямая из этих двух пучков в точности соответствует точке пересечения
  этой прямой с прямой $\langle p_2, p_3 \rangle$, откуда получаем, что параметризация $N(\alpha, \beta)$ естественным образом заменяется на параметризацию $N(s_0, s_1)$ упорядоченными парами точек
  на прямой $\langle p_2, p_3 \rangle$. Наконец, если прямая $l \subset \mathbb{C} \mathbb{P}^3$ одновременно пересекает прямые $\langle p_0, s_0 \rangle$ и $\langle p_1, s_1 \rangle$, то отсюда следует что она является
  пересечением  плоскостей $\langle l, p_0 \rangle$ и $\langle l, p_1 \rangle$, откуда видно, что пучки $Q_0$ и $Q_1$ могут быть заменены на соответствующие пучки плоскостей, проходящих через
  $\langle p_i, s_i \rangle$, следовательно аналитическое описание, сопоставляющее точке $[l] \in {\rm Gr}(1, 3)$ упорядоченную пару точек $(s_0, s_1)$ на прямой $\langle p_2, p_3 \rangle$
  указанным выше способом в точности соответствует геометрическому описанию ${\rm Gr}^0(1, 3)$, представленному выше.

  Заметим далее, что для корректности отождествления в аналитической ситуации ${\rm Gr}^0(1, 3)$ с расслоением нам необходимо исключить базисные множества пучков $Q_0$ и $Q_1$,
   которые представлены в виде $B_0 = \{ w_{02} = w_{03} = 0 \}$ и $B_1 = \{w_{12} = w_{13} = 0 \}$. Нетрудно видеть, что $B_0$ эквивалентно условию $l \subset \langle p_0, p_2, p_3 \rangle$,
    а $B_1$ эквивалентно условию $l \subset \langle p_1, p_2, p_3 \rangle$, и оба эти случая исключены из $Gr^0(1, 3)$ выбором некритических значений отображений моментов.
     
     Пусть  $N^0(s_0, s_1) = N(s_0, s_1) \cap {\rm Gr}^0(1,3)$,
    тогда  очевидным образом имеем расслоение $\tau: {\rm Gr}^0(1,3) \to \langle p_2, p_3 \rangle \times \langle p_2, p_3 \rangle$ со слоем $N^0(s_0, s_1)$. Аналитически эти слои представляются
    как пересечения  $Q_0 \cap Q_1 \cap {\rm Gr}^0(1, 3)$, то есть задаются тремя уравнениями
    $$
    \alpha_0 w_{02} + \alpha_1 w_{03} = 0, \quad \beta_0 w_{12} + \beta_1 w_{13} = 0, \quad  w_{01}w_{23} - w_{02} w_{13} + w_{03} w_{12}=0,
    $$
    откуда видно, что если точки $s_0, s_1 \in  \langle p_2, p_3 \rangle$ не совпадают, то пересечение представляет собой невырожденную 2 - мерную квадрику, если же
    $s_0 = s_1$, что эквивалентно $\alpha_0 \beta_1 = \alpha_1 \beta_0$, то пересечение представляется парой плоскостей $\{ w_{01} = 0 \} \cup \{w_{23} = 0 \}$.
    Заметим при этом, что  это пересечение всегда является торическим, поскольку имеется два отображения моментов $\tilde \mu_0, \tilde \mu_1$,
    сохраняющих каждый слой гамильтоновым действием. Из формулы (2) получаем
    $$
    \tilde \mu_0 + \tilde \mu_1 = 1 +  \frac{\vert w_{01} \vert^2 - \vert w_{23} \vert^2}{\sum_{i<j} \vert w_{ij} \vert^2},
    $$
    откуда видно, что $\tilde \mu_0 + \tilde \mu_1 < 1$ соответствует выбору плоскости $ w_{01} = 0$ в вырожденном слое над диагональю в $\langle p_2, p_3 \rangle \times \langle p_2, p_3 \rangle$.

    Таким образом, глобально картина выглядит так: над общей точкой прямого произведения $\langle p_2, p_3 \rangle \times \langle p_2, p_3 \rangle$ имеем торический слой --- невырожденную квадрику,
    многоугольником Дельцана которой является квадрат; над точкой, лежащей на диагонали, торический слой вырождается в объединение пары плоскостей,
    многоугольниками Дельцана которой будут треугольники, на которые распадается квадрат; пересечению этих двух плоскостей соответствует диагональ квадрата
    $\tilde \mu_0 + \tilde \mu_1 = 1$. Отсюда видно, что $N^0(s_0, s_1)$ есть прообраз открытого треугольника так что прообразом каждой точки этого открытого
    треугольника является гладкий двумерный тор, инвариантный относительно гамильтонова действия ${\rm T}^2$, натянутого на коммутирующие функции $\tilde \mu_0, \tilde \mu_1$.
    Этим завершается доказательство Предложения 1.

    Нетрудно видеть, что утверждение Теоремы 2  для случая ${\rm Gr}(1, 3)$  немедленно следует из Предложения 1.

\section{Общий случай}

Рассмотрим теперь общий случай грассманиана ${\rm Gr}(1, n)$ прямых в $\mathbb{C} \mathbb{P}^n$. Зафиксируем однородные координаты $[z_0: ... : z_n]$, согласованные с кэлеровой структурой,
и рассмотрим $k$ стандартных  отображений моментов
$$
\mu_i = \frac{\vert z_i \vert^2}{\sum_{i=0}^n \vert z_i \vert^2}, \quad i = 0, ..., k-1,
$$
 и соответствующие им отображения моментов $\tilde \mu_i$ на ${\rm Gr}(1, n)$. Снова выделим в ${\rm Gr}(1, n)$ открытую компоненту ${\rm Gr}^0(1, n)$, накладывая условия $\tilde \mu_i > 0, \tilde \mu_0 + ... +  \tilde \mu_{k-1} < 1$ и покажем, что

 {\bf Предложение 2.} {\it Открытое многообразие ${\rm Gr}^0(1, n)$ расслоено над многообразием $M_{n-k}$, являющимся тотальным пространством $k$ - кратного прямого произведения проективизации  тавтологического расслоения $\mathbb{P}(\tau) \to {\rm Gr}(1, \langle p_k, ..., p_{n+1} \rangle )$ на себя, так что гамильтоново действие тора $T^k$, натянутое на
  $\tilde \mu_0, ..., \tilde \mu_{k-1}$, сохраняет базу и слои.}

 Здесь и ниже как  ${\rm Gr} (1, H)$ мы обозначаем многообразие проективных прямых в проективном пространстве $H$, в том числе когда $H$ есть проективная оболочка точек $\langle p_{i_1}, ... p_{i_m} \rangle$.

 {\bf Доказательство.} Так как мы снова, как и в случае ${\rm Gr}(1, 3)$ выше, исключили из рассмотрений те прямые $l \subset \mathbb{C}\mathbb{P}^n$, которые пересекают подпространство $\langle p_0, ...,  p_{k-1} \rangle $  или лежат в гиперплоскостях 
 $$
 \langle p_1, p_2, ..., p_n \rangle, \langle p_0, p_2, ..., p_n \rangle, ...., \langle p_0, ..., p_{k-2}, p_k, ..., p_n \rangle,
 $$
  то ${\rm Gr}^0(1, n)$ можно расслоить над
 $M_{n-k}$: любая прямая $l$ из ${\rm Gr}^0(1, n)$ скрещивается с $\langle p_0, ..., p_{k-1} \rangle$,  поэтому проективная оболочка $H_l = \langle l, p_0, ...,  p_{k-1} \rangle$ обязательно проективно  $k+1$ - мерна; с другой стороны, $l$ однозначно определяет свое подпространство $H_l$ и прямую $l_0 \subset H_l \cap \langle p_k, ..., p_n \rangle$ в последнем $(n-k)$ - мерном проективном пространстве. Нетрудно видеть, что само $H_l$ однозначно восстанавливается по последней прямой $l_0$, так как является проективной оболочкой
   $\langle l_0, p_0, ..., p_{k-1} \rangle$.

  Тогда имеется следующее замечание:  гамильтоново действие ${\rm T}^k$, порождаемое выбранными нами  отображениями моментов $\tilde \mu_0, ..., \tilde \mu_{k-1}$, сохраняет подмножества ${\rm Gr}^0(1, H_l) \subset {\rm Gr}^0(1, n)$,   где ${\rm Gr}^0(1, H_l)$ есть множество прямых в $n-k$ - мерном проективном подпространстве $H_l  \subset \mathbb{C} \mathbb{P}^n$, удовлетворяющих условиям на функции $\tilde \mu_i$,
  выделяющих компоненту ${\rm Gr}^0(1, n)$. В самом деле, $k+1$ - мерное подпространство $H_l$
  однозначно определяется парой проективных подпространств $l_0, \langle p_0, ..., p_{k-1} \rangle$, лежащих в нем. Гамильтоново действие $\tilde \mu_i$ оставляет прямую $l_0 \subset \langle p_k, ..., p_n \rangle$
  неподвижной поскольку она полностью содержится в критическом множестве, соответствующем нулевым критическим значениям всех рассматриваемых  отображений моментов $\tilde \mu_0, ...,
  \tilde \mu_{k-1}$; с другой стороны,
  проективное подпространство  $\langle p_0, ... , p_{k-1} \rangle$ также инвариантно (хотя и не неподвижно), откуда следует наше замечание.

  Таким образом, наша задача расщепляется на две: первая --- симплектическая редукция по гамильтонову действию на каждом грассманиане ${\rm Gr}(1, H_l)$ при фиксированном $H_l$, инвариантном относительно действия тора $T^k$, натянутого на $\tilde \mu_0, ..., \tilde \mu_{k-1}$,    вторая --- глобализация конструкции и описание подходящего параметризующего многообразия.

  Первая задача  решается по образцу из первого параграфа: все множество ${\rm Gr}(1, H_l)$ с набором отображений  $\tilde \mu_i|_{{\rm Gr}(1, H_l)}, i = 0, ..., k-1$,
   естественным образом расслаивается на торические слои. В самом деле, сопоставим прямой $l \subset H_l$ набор точек $(s_0, ..., s_{k-1}), s_i \in l_0$, по правилу:
  сначала рассмотрим точку  $s_j^0 = \langle l_0, p_{m_1}, ..., p_{m_{k-1}} \rangle \cap l$, где $(j, m_1, ..., m_{k-1})$ есть перестановка чисел $(0, ..., k-1)$, на прямой $l$, и затем спроектируем ее
  из $\langle p_{m_1}, ..., p_{m_{k-1}} \rangle$ на $l_0$ в проективном пространстве $\langle l_0, p_{m_1}, ..., p_{m_{k-1}} \rangle$, обозначая результат как $s_j \in l_0$. Тогда в ${\rm Gr}(1, H_l)$ имеем набор дивизоров
   $\{ D(s_j) \}$, определяемых условием
   $$
   D(s_j) = \{ [l] \in {\rm Gr}(1, H_l) \quad  | \quad  l \cap \langle s_j, p_{m_1}, ..., p_{m_{k-1}} \rangle \neq \emptyset \}
   $$
   (стандартный способ определения дивизоров в грассманиане, см. [4]).

   Заметим, что для любой точки $s_j \in l_0$ соответствующий дивизор $D(s_j)$ инвариантен относительно гамильтонова действия каждого нашего отображения моментов $\tilde \mu_i$.
   В самом деле, наше  проективное подпространство $\langle s_j, p_{m_1}, ..., p_{m_{k-1}} \rangle \subset H_l$ натянуто на точки, неподвижные относительно гамильтонова действия $\mu_i$, отсюда любая прямая
   $[l] \in D(s_j)$ под дейcтвием $\tilde \mu_i$ переходит в прямую из того же подмногообразия. Отсюда следует, что подмногообразие
   $$
   Y(s_0, ..., s_{k-1}) = \bigcap_{j =0}^{k-1} D(s_j) \subset {\rm Gr}(1, H_l)
   $$
    инвариантно относительно действия всего тора $T^k$, натянутого на отображения моментов $\tilde \mu_0, ... , \tilde \mu_{k-1}$; с другой стороны, размерность этого подмногообразия равна $k$, откуда следует, что для каждого набора $(s_0, ..., s_{k-1})$ имеем торический слой
   $Y(s_0, ..., s_{k-1})$ (не обязательно гладкий).

   Согласно общей теории торических многообразий, см. [5], наше многообразие $Y(s_0, ..., s_{k-1})$, снабженное набором отображений моментов $\tilde \mu_0, ..., \tilde \mu_{k-1}$,
   однозначно определяется своим выпуклым многогранником Дельцана $P_Y \subset \mathbb{R}^k$, являющимся образом отображения действия
   $$
   F_{act} = (\tilde \mu_0, ... \tilde \mu_{k-1}): Y(s_0, ..., s_{k-1}) \to \mathbb{R}^k.
   $$

   Нетрудно видеть, что для общего набора $(s_0, ..., s_{k-1})$ многогранник Дельцана $P_Y$  задается так: рассмотрим единичный куб в $\mathbb{R}^k$ и отрежем от него гиперплоскостью
   $\sum_{i=1}^k x_i = 2$ ту часть, которая содержит начало координат. В самом деле, по  определению торическое многообразие $Y(s_0, ..., s_{k-1})$ содержит в себе прямые,
   являющиеся прообразами вершин вида (все $x_i$ нули), (ровно одно $x_i$ отлично от нуля), (ровно два $x_i, x_j$ отличны от нуля), но именно для такого набора вершин единичного куба выпуклой оболочкой будет $P_Y$, представленный выше. Предьявим такие прямые: прямая $l = l_0 \subset H_l$ переходит в точку первого вида; прямая $\langle p_i, s_i \rangle \subset H_l, i = 0, ..., k-1$, переходит в точку, у которой $i$ - ая координата равна 1, а все остальные нули; наконец, прямая $\langle p_i, p_j \rangle \subset H_l, 0 \leq i < j \leq k-1$, переходит в вершину, у которой $i$ - ая и $j$ - ая координаты равны 1,
а все остальные нули.

Как и выше в случае ${\rm Gr}(1, 3)$, при совпадении пары точек $s_i = s_j$, определяющих многообразие $Y(s_0, ..., s_{k-1})$, это торическое многообразие становится приводимым,
состоящим из торических компонент. В самом деле, в этом случае пересечение дивизоров $D(s_i) \cap D(s_j)$  распадается в объединение пары компонент $(D(s_i) \cap D(s_j))_a \cup
(D(s_i) \cap D(s_j))_b$, где $a$ компонента состоит из прямых, проходящих через точку $s_i = s_j$, а $b$ - компонента состоит из прямых, содержащихся в проективном подпространстве
  $ \langle p_0, ..., p_{k-1}, s_i= s_j \rangle \subset H_l$. Как и в случае ${\rm Gr}(1, 3)$ эти компоненты различаются значениями отображений моментов: $a$ - компонента соответствует
  случаю $\sum_{i=0}^{k-1} \tilde \mu_i \leq 1$, а $b$ - компонента соответствуют случаю, когда та же сумма больше или равна 1.

  При совпадении большего числа $s_i$ мы получаем более специальную картину распадения $Y(s_0, ..., s_{k-1})$ на компоненты, однако полное исследование соответствующей комбинаторной структуры  не входит в круг задач настоящей работы, поскольку для доказательства Теоремы 2 нам достаточно восстановить наши условия $c_i > 0, \sum_{i=0}^{k-1} c_i < 1$,
  накладываемые на значения отображений моментов $\tilde \mu_0, ..., \tilde \mu_{k-1}$ и убедиться в том, что каждый такой набор $(c_0, ..., c_{k-1})$
  определяет в {\bf каждом } торическом многообразии $Y(s_0, ..., s_{k-1})$ соответствующий гладкий $k$ - мерный тор Лиувилля, на котором $T^k$, порождаемый
  гамильтоновым действием $\tilde \mu_0, ..., \tilde \mu_{k-1}$, действует свободно.

  Таким образом в качестве промежуточного вывода мы получаем, что результатом симплектической редукции ${\rm Gr}(1, H_l)= {\rm Gr}(1, k+1)$ по действию $k$ - мерного
  тора при выборе значений отображений моментов  $c_i > 0, \sum_{i=0}^{k-1} c_i < 1$ есть  $k$ - кратное прямое произведение проективной прямой $l_0$ на себя.

   Перейдем теперь ко второму шагу: очевидно, что сопоставление $l \mapsto H_l \mapsto l_0 \subset \langle p_k, ..., p_n \rangle$ взаимно - однозначно по второй стрелке, поэтому
   отсюда получаем отображение ${\rm Gr}^0(1, n) \to {\rm Gr}(1, \langle p_k, ..., p_n \rangle)$. Как было показано выше, над каждой точкой $[l_0] \in {\rm Gr} (1, \langle p_k, ..., p_n \rangle)$    имеем расслоение на торические открытые  части $Y^0(s_0, ..., s_{k-1}) \subset Y(s_0, ..., s_{k-1})$, параметризованные набором точек $(s_0, ..., s_{n-1})$,
   каждая из которых лежит на $l_0$. Отсюда, глобализуя над всем ${\rm Gr}(1, \langle p_k, ..., p_n \rangle)$ эту картину, с учетом того, что точки $s_i$ оказываются
   точками на проективизации тавтологического расслоения $\tau \to {\rm Gr}(1, n-k)$, получаем утверждение Предложения 2. Мы будем обозначать соответствующее
   многообразие ${\rm tot}(\mathbb{P}(\tau) \times ... \mathbb{P}(\tau) \to {\rm Gr}(1, n-k))$ как $M_{n-k}$.

    Если теперь мы выберем набор значений $c_0, ..., c_{k-1} \in \mathbb{R}_+$ так
   что сумма $\sum_{i=0}^{k-1} c_i  < 1$, то этим в каждом торическом слое ${\rm Gr}^0(1, n) \to M_{n-k}$ будет выделен гладкий тор, который при факторизации даст точку,
   откуда следует утверждение Теоремы 2:
   $$
   {\rm Gr}(1, n) \slash \slash {\rm T}^k \cong {\rm tot} (\mathbb{P}(\tau) \times ... \times \mathbb{P}(\tau) \to {\rm Gr}(1, n-k)),
   \eqno (3)
   $$
   где в слое $k$ - кратное прямое произведение.

   {\bf Замечание 2.} Приведенное выше представление не позволяет ответить на вопрос о том, какая симплектическая форма получается в процессе редукции из исходной кэлеровой формы на ${\rm Gr}(1, n)$ на редуцированном   симплектическом многообразии $M_{n-k}$. Несмотря на то, что при близких значениях отображений моментов мы получаем изоморфные многообразия в результате симплектической редукции, симплектические формы на таких многообразиях отличаются друг от друга. Например, мы не можем утверждать в общем случае, что ограничение симплектической формы $\omega_t, t = (c_0, ..., c_{k-1})$, получаемой в результате симплектической редукции на $M_{n-k}$, на проективные прямые $\mathbb{P}(\tau)$ над одной и той же точкой базы ${\rm Gr}(1, n-k)$, дают одну и ту же симплектическую форму. Для части построенных подмногообразий нам это и не потребуется. Однако поскольку ниже мы используем нестандартные антиголоморфные инволюции, то нам необходимо следующее простое замечание. Так как при автоморфизме исходного проективного пространства $\mathbb{C} \mathbb{P}^n$, индуцируемого перестановкой однородных координат $z_i \leftrightarrow z_j$, имеем: 1) симметрию отображений моментов $\tilde \mu_i \leftrightarrow \tilde \mu_j$; 2) симметрию точек $s_i \leftrightarrow s_j$ на проективной прямой $l_0$; 3) а значит и кэлерову симметрию между $i$ -ыми и $j$ -тыми прямыми слагаемыми в $M_{n-k}$ над базой ${\rm Gr}(1, n-k)$. Отсюда видно, что если в дальнейшем мы дополним условия на значения отображений моментов дополнительным условием симметрии $c_0 = ... = c_{k-1}$, то этим на $M_{n-k}$ будет определено действие группы перестановок кэлеровыми изометриями.

   Однако заметим при этом, что для реализации конструкции Миронова, где важную роль имеют вещественные части многообразий с голоморфной инволюцией, нам не обязательно рассматривать 
   полностью симметричный случай $c_0 = ... = c_{k-1}$. В самом деле,  наши отображения моментов $\tilde \mu_i$
   действуют кэлеровыми изометриями на ${\rm Gr} (1, n)$, а группа ${\rm T}^k$ коммутативна, поэтому любое значение отображение моментов (в смысле ${\frak g}^*$) лежит в центре группы, что требуется для    корректной реализации кэлеровой редукции, см. [4]. Так что построенное нами многообразие $M_{n-k}$ самой конструкцией наделено кэлеровой структурой, и при этом такая структура согласована с антиголоморфной инволюцией $\tau: M_{n-k} \to M_{n-k}$, также индуцируемой из стандартной антиголоморфной инволюции на $\mathbb{C} \mathbb{P}^n$. В самом деле,
    нетрудно видеть, что стандартное сопряжение прямой $l \subset \mathbb{C} \mathbb{P}$ согласовано с сопряжением тех данных, которым мы эту прямую сопоставили: так как все точки $p_i 
    \in \mathbb{C} \mathbb{P}^n$ вещественны, то $\bar l$ будет соответствовать $\bar l_0 \subset \langle p_k, ..., p_n \rangle$, и набор точек $s_i \in l_0$ также преобразуется в сопряженный набор $\bar s_i \in \bar l_0$. Поэтому какова бы ни была кэлерова форма на $M_{n-k}$ в результате кэлеровой редукции, вещественная часть этого многообразия будет
    лагранжевым подмногообразием; и эта вещественная часть выделена условиями вещественности, налагаемыми и на базу ${\rm Gr}(1, n-k)$ (чтобы $l_0$ была вещественной),
    и на точки в слоях всех прямых слагаемых $\mathbb{P} (\tau)$. Иными словами наше  фактор - многообразие $M_{n-k}$ естественно вкладывается в прямое произведение
   $$
   {\rm Gr}(1, \langle p_k, ..., p_n \rangle) \times \langle p_k, ..., p_n \rangle \times ... \times \langle p_k, ..., p_n \rangle,
   \eqno (4)
   $$
   где справа имеем прямое произведение $k$ копий одинаковых проективных пространств с одинаковыми комплексными структурами, порождаемыми исходной комплексной  структурой
   на $\langle p_k, ..., p_n \rangle$, с одинковыми антиголоморфными инволюциями, но каждая копия со своей кэлеровой (симплектической) формой. Редуцированное  многообразие $M_{n-k}$ реализуется как пересечение $k$ циклов
   $$
   M_{n-k} = \bigcap_{i=1}^k \tilde {\cal U}_i,
   \eqno (5)
   $$
    где $\tilde {\cal U}_i  = \langle p_k, ..., p_n \rangle \times ... \times {\cal U}_i \times ... \times  \langle p_k, ..., p_n \rangle$, получаемое подстановкой в прямое произведение на $i$ - ом месте
     ${\cal U}_i \subset  {\rm Gr}(1, \langle p_2, ..., p_n \rangle) \times \langle p_k, ..., p_n \rangle$ --- цикла инциденции для $i$ - го проективного пространства;
    и редуцированная кэлерова структура есть ограничение  с прямого произведения на подмногообразие $M_{n-k}$, определяемое формулой (4).

   {\bf Замечание 3.} Нетрудно видеть, что как и в случае ${\rm Gr}(1,3)$, расслоение ${\rm Gr}^0(1, n) \to M_{n-k}$ топологически не тривиально.
   В самом деле, вернемся к дивизорам $D(s_i) \subset Gr(1, H_l)$ для набора точек $s_0, ..., s_{k-1} \in l_0$. Вместо $D(s_i)$ можно рассмотреть
   пучок $k$ - мерных проективных подпространств, содержащих $k-1$ - мерное подпространство $ \langle s_j, p_{m_1}, ..., p_{m_{k-1}} \rangle$: очевидно, что если
   прямая $l$ пересекает $\langle s_j, p_{m_1}, ..., p_{m_{k-1}} \rangle$, то она обязана лежать в одном из таких $k$ - мерных подпространств.
   Отсюда имеем расслоение над $k$ кратным прямым произведением $l_0$ на себя, слоем которого будет $k$ - кратное прямое произведение проективных прямых,
   параметризующих описанные пучки. Однако каждый такой пучок имеет такой же вид, как и в случае ${\rm Gr}(1, 3)$, см. Замечание  1 в предыдущем параграфе.
   В самом деле, на проективной плоскости $\langle p_j, l_0 \rangle$ пучок соответствующих $k$ - мерных проективных подпространств высекает пучок прямых, проходящих через $s_j$,
   и варьируя $s_j$ вдоль $l_0$ снова получаем несколько копий расслоений Хопфа в случае фиксации подходящего набора значений $c_i> 0, \sum_{i=0}^{k-1} c_i < 1$
   отображений моментов $\tilde \mu_i$.  Очевидно, что конструкция глобализуется и над
   ${\rm Gr}(1, n-k)$, откуда следует топологический тип расслоения ${\rm Gr}^0(1, n) \to M_{n-k}$.

   {\bf Замечание 4.} В некоторых задачах (например, в главной задаче работы [1] -- построения минимальных и  гамильтоново минимальных лагранжевых подмногообразий) необходимо бывает рассмотрение
   конкретных выделенных значений $t \in {\frak g}^*$ отображения моментов, для которых на соответствующем многообразии уровня $N_t$ группа действует не свободно,
   что приводит к возникновению особенностей соответствующего фактор - пространства $N_t \slash G$. Например, в наших рассмотрениях Предложения 1 можно было
   допустить выбор значений $c_0, c_1$ отображений моментов $\tilde \mu_0, \tilde \mu_1$ так что $c_0 + c_1 = 1$. В этом случае, как следует из формул параграфа 1,
   над диагональю в прямом произведении $\langle p_2, p_3 \rangle \times \langle p_2, p_3 \rangle$ вместо двумерных торов будут лежать одномерные торы, то есть окружности, лежащие на прямых из пересечения $\{w_{01} = 0 \} \cap \{ w_{23} = 0 \}$. Поэтому при факторизации для таких значений отображений моментов мы получим редуцированное многообразие с особенностями,
   однако конструкция Теоремы 1 может быть применена и в этом случае, если рассматривать лагранжевы погружения, в том числе и лагранжевы погружения в симплектическое
   многообразие с особенностями. Эта отдельная тема, вытекающая из Теоремы 1, требует и отдельного внимательного исследования.

   {\bf Замечание 5.} Наша конструкция не может быть прямо перенесена на случаи $k=1$ и $k=n$: в первом случае такой необходимости нет, поскольку в работе [3] мы уже построили
   соответствующее лагранжево подмногообразие, которое оказалось также того же вида --- расслоение на двумерные торы над ${\rm Gr}(1, n-1)$; в случае $k=n$ дополнительное
   к $\langle p_0, ..., p_{n-1} \rangle$ проективное пространство нульмерно, поэтому прямые некуда проектировать. Поскольку в настоящей работе мы пользумеся именно представленным выше методом,
   отложим рассмотрения случая $k = n$ до будущей работы.

   Теперь мы готовы применить Теоремы 1 и 2 для построения большого запаса примеров лагранжевых подмногообразий в грассманиане ${\rm Gr}(1, n)$, что является нашей главной задачей.

\section{Примеры}

 {\bf Предложение 3.} {\it Многообразие ${\rm Gr} (1, n)$, снабженное симплектической формой плюккеровым вложением, допускает набор гладких лагранжевых подмногообразий
 $\{ S_k \}, k =0, ..., n-1$, представляемых как гладкие расслоения на $T^{2k}$ - мерные торы над вещественными грассманианами ${\rm Gr}_{\mathbb{R}} (1, n-k)$.}

{\bf Доказательство.} Начнем со случаев $k = 0$ и  1: первый очевиден, так как утверждение в этом случае есть просто тавтология; второй доказан в работе [3].
Далее, из Теоремы 2 следует, что для каждого $k = 2, ..., n-1$, результатом симплектической редукции ${\rm Gr}(1, n) \slash \slash T^k$ при указанном выше
выборе отображений моментов и их значений является кэлерово многообразие $M_{n-k}$. Вещественная часть $M^{\mathbb{R}}_{n-k}$ относительно антиголоморфного сопряжения,
индуцируемого стандартными антиголомофрными сопряжениями на $\mathbb{C} \mathbb{P}^k$ и ${\rm Gr}(1, n-k)$ в прямом произведении (4), представляется как вещественная версия
пересечения (5), откуда получаем
$$
M^{\mathbb{R}}_{n-k} = {\rm tot} (\mathbb{P}(\tau_{\mathbb{R}}) \times ... \times \mathbb{P} (\tau_{\mathbb{R}}) \to {\rm Gr}_{\mathbb{R}}(1, n-k)),
$$
и поскольку слой $\mathbb{P}(\tau_{\mathbb{R}}) \to {\rm Gr}_{\mathbb{R}}(1, n-k)$ есть окружность, то глобально $M^{\mathbb{R}}_{n-k}$ есть топологически нетривиальное
гладкое расслоение со слоем $S^1 \times ... \times S^1 = T^k$ (исключение составляет случай $k = n-1$, когда ${\rm Gr}(1, n-k)$ есть просто точка,
и $M^{\mathbb{R}}_{n-1}$ просто изоморфно $T^{n-1}$).

Из Теоремы 1  и Замечания 3 мы знаем, что каждому лагранжеву подмногообразию $L \subset M_{n-k}$ естественно соответствует лагранжево подмногообразие в ${\rm Gr}(1, n)$,
топологический тип которого есть $T^k$ расслоение над  $L$: поскольку $M^{\mathbb{R}}_{n-k} \subset M_{n-k}$ лагранжево, то для доказательства Предложения 3 нам достаточно показать, что несмотря на топологическую нетривиальность $T^k$ над $M_{n-k}$,  соответствующее ему лагранжево подмногообразие имеет тип расслоения над ${\rm Gr}_{\mathbb{R}}(1, n-k)$
на $2k$ - мерные торы. Однако как было показано в Замечании 3 слой над каждой точкой ${\rm Gr}_{\mathbb{R}}(1, n-k)$ представляется как $k$ - копий расслоений Хопфа,
которые при ограничении на вещественный тор $\mathbb{P} (\tau_{\mathbb{R}}) \times ... \mathbb{P}(\tau_{\mathbb{R}})$, как указано в Замечании 1, дают
прямое произведение $k$ копий $T^2$, то есть тор $T^{2k}$. Подчеркнем, что по построению соответствующее расслоение $S_k \to {\rm Gr}_{\mathbb{R}}(1, n-k)$
на $2k$ - мерные торы топологически нетривиально, аналогично тому, что было получено для случае $k=1$ в работе [3].  Этим завершается доказательство Предложения 3.

{\bf Замечание 6.} Нетрудно видеть, что представленные выше лагранжевы подмногообразия $S_k \subset {\rm Gr}(1, n)$ есть в точности те, которые
дает обобщенная конструкция Миронова в случае отображений моментов $\tilde \mu_i$ и выбираемых выше значениях для них. Поэтому в настоящей работе мы
в достаточной степени продвинулись в задаче описания типов лагранжевых подмногообразий в ${\rm Gr} (1, n)$ с помощью обобщенной конструкции Миронова, сформулированной в работе [2].
Кроме того, все лагранжевы подмногообразия $\{S_k \}$ связаны между собой преобразованиями, называемыми в Замечании 1 лагранжевыми перестройками, а именно имеем цепочку
$$
{\rm Gr}_{\mathbb{R}}(1, n) = S_0 \mapsto S_1 \mapsto ... \mapsto S_{n-1},
$$
где $i$ - ая стрелка соответствует лагранжевой перестройке, индуцируемой отображением моментов $\tilde \mu_{i-1}$; одновременно и для любой пары $0 \leq i < j \leq n-1$
соответствующий набор отображений моментов $\tilde \mu_i, ..., \tilde \mu_{j-1}$ индуцирует лагранжеву перестройку $S_i \mapsto S_j$.

Предложение 3 превращает Гипотезу из работы [2] в утверждение, однако формулы (4) и (5) подсказывают еще целый класс многообразий, существенно расширяющих
множество примеров лагранжевых подмногообразий.

Начнем как всегда с исследования наиболее простой ситуации и рассмотрим случай $k=2$. Тогда, как в (4) и (5),  реализуем $M_{n-2}$ как подмногообразие в тройном прямом произведении $\mathbb{C} \mathbb{P}^{n-2} \times {\rm Gr}(1, n-2) \times \mathbb{C} \mathbb{P}^{n-2}$, где первое и третье прямое слагаемое отождествлены как комплексные многообразия; если ввиду Замечания 2 дополнить наши условия на значения отображений моментов еще одним, а именно  $c_0 = c_1$, то первое и второе слагаемое отождествлены и как кэлеровы многообразия.
Заметим, что для этого случая многообразие $M_{n-2}$ бирационально изоморфно прямому произведению $\mathbb{C} \mathbb{P}^{n-2}  \times \mathbb{C} \mathbb{P}^{n-2}$,
будучи раздутием диагонали в этом прямом произведении. В последнем прямом произведении известен пример лагранжева подмногообразия, представляющего нетривиальный класс гомологий, а именно
 введя однородные координаты $[x_i],  [y_i]$ на первом и втором прямых слагаемых,
рассмотрим  лагранжево подмногообразие --- антидиагональ, задаваемую условием
$$
\bar \Delta = \{ y_i = \bar x_i \} \subset \mathbb{C} \mathbb{P}^{n-2} \times \mathbb{C} \mathbb{P}^{n-2}.
$$
Поскольку $\bar \Delta$ пересекается с $\Delta$  в вещественных точках, то полный прообраз $\bar \Delta$ в $M_{n-2}$ относительно проекции
$$
\pi: M_{n-2} \to \mathbb{C} \mathbb{P}^{n-2} \times \mathbb{C} \mathbb{P}^{n-2}
$$
очевидно не может быть лагранжевым в $M_{n-2}$: слои $\pi$ над вещественными точками являются комплексными подмногообразиями.
Однако используя представления (4) и (5) мы можем построить гладкое лагранжево подмногообразие в $M_{n-2}$ аналогичного типа.
А именно, на прямом произведении (4) рассмотрим антиголоморфную инволюцию, действующую стандартно на ${\rm Gr} (1, n-2)$ и
антидиагональным образом на $\mathbb{C} \mathbb{P}^{n-2} \times \mathbb{C} \mathbb{P}^{n-2}$ то есть так что
$([x_i], [y_i]) \mapsto ([\bar y_i], [\bar x_i])$. Очевидно, что такая антиголоморфная инволюция в качестве неподвижных точек имеет
лагранжево подмногообразие в (4)(здесь важным является совпадение кэлеровых форм, что выполнено при $c_0 = c_1$), однако замечательным образом она естественно ограничивается и на $M_{n-2}$ и задает нестандартную антиголоморфную инволюцию и здесь. Покажем, что $M_{n-2}$ инвариантно относительно такой инволюции откуда и будет следовать необходимый факт.
В самом деле, по  определению точкой $M_{n-2}$ является набор $(p_1, l, p_2)$ из прямой $l$ и двух точек на ней; при нашем нестандартном сопряжении
этот набор переходит в набор $(\bar p_2, \bar l, \bar p_1)$, который очевидно обладает тем же свойством инциденции, и, следовательно,
нестандартная инволюция сохраняет $M_{n-2}$. Следовательно, если множество неподвижных точек $M^r_{n-2} \subset M_{n-2}$ имеет правильную
размерность, то $M^r_{n-2}$ является лагранжевым подмногообразием. Однако если вспомнить представление
$$
M_{n-2} \equiv {\rm tot}(\mathbb{P} (\tau) \times \mathbb{P}(\tau) \to {\rm Gr}(1, n-2))
$$
как расслоения над ${\rm Gr}(1, n-2)$, то нетрудно видеть, что $M^r_{n-2}$ есть расслоение над ${\rm Gr}_{\mathbb{R}}(1, n-2)$
со слоем $S^2$. Вещественная размерность этого тотального пространства в точности равна комплексной размерности $M_{n-2}$, откуда получаем
новый пример лагранжева подмногообразия в последнем, а значит и в ${\rm Gr} (1, n)$.

Определим топологический тип такого лагранжева подмногообразия. Поскольку для фиксированных подходящих значений $\tilde \mu_0, \tilde \mu_1$
мы имеем в качестве слоя над точкой ${\rm Gr}(1, n-2)$ две копии расслоения Хопфа, то необходимо установить топологический тип ограничения расслоения Хопфа на
антидиагональ в прямом произведении $\mathbb{P}(\tau) \times \mathbb{P} (\tau)$. Однако такое расслоение должно иметь нулевой класс Черна, а $S^2$ односвязна,
откуда следует, что над каждой точкой ${\rm Gr}_{\mathbb{R}}(1, n-2)$ имеется слой $T^2 \times S^2$, а значит глобально мы получаем
расслоение над ${\rm Gr}_{\mathbb{R}}(1, n-2)$ со слоем  $T^2 \times S^2$. При этом очевидно, что такой топологический тип отличается от топологических типов
 $S_k$ из Предложения 3  для любого $i$. Следовательно, список топологических типов гладких лагранжевых подмногообразий
пополнен еще одним типом, возникающим при рассмотрениях редукции по $T^2$.

Но почему бы нам не реализовать такую же возможность в случае $k = 3$? Ведь в этом случае $M_{n-3}$ изоморфно
${\rm tot} (\mathbb{P}(\tau) \times \mathbb{P}(\tau) \times \mathbb{P}(\tau) \to {\rm Gr}(1, n-3))$, и поэтому
допускает пример подмногообразия смешанного типа, а именно над вещественной частью грассманиана ${\rm Gr}_{\mathbb{R}}(1, n-3)$
рассмотрим из одного прямого слагаемого вещественную часть $\mathbb{P}(\tau_{\mathbb{R}})$, а с двумя другими компонентами
проделаем то же, что и в случае $k=2$. Очевидно, что в результате мы получим еще одно лагранжево подмногообразие
в $M_{n-3}$, тип которого есть топологически нетривиальное расслоение со слоем $S^2 \times S^1$ над ${\rm Gr}_{\mathbb{R}}(1, n-3)$.
Далее, обобщая на случай $k = 4$ получаем еще больше возможностей: над ${\rm Gr}_{\mathbb{R}}(1, n-4)$ имеем расслоения
со слоем $S^1 \times S^1 \times S^1 \times S^1$ (обобщенный цикл Миронова $S_4$ из Предложения 3), со слоем $S^2 \times S^1 \times S^1$
и со слоем $S^2 \times S^2$, то есть 3 разных топологических типа. И естественно, что каждый такой тип определяет после подъема вдоль соответствующего $T^k$
- расслоения соответствующий тип лагранжева подмногообразия в ${\rm Gr}(1, n)$. Единственное требование, необходимое для данных построений, заключается в совпадении
соответствующих значений отображений моментов. Но если просто положить $1> c_0 = ... = c_{k-1} > 0, c_0 < 1/k$, то этим достигается возможность рассматривать любые пары. 

Из представленных конструкций и рассуждений  получаем

{\bf Следствие.} {\it Число различных топологических типов лагранжевых подмногообразий в грассманиане ${\rm Gr}(1, n)$ не меньше $n + [\frac{n}{2}] \cdot [\frac{n-1}{2}]$.}

{\bf Доказательство.} Для каждого фиксированного $k = 0, ..., n-1$, согласно нашим конструкциям и замечаниям, мы имеем не менее $1 + [\frac{k}{2}]$
разных типов лагранжевых подмногообразий в ${\rm Gr}(1, n)$ вследствии Теоремы 1. В самом деле, для $M_{n-k}$ имеется стандартная вещественная часть
$M^{\mathbb{R}}_{n-k}$ ( что дает единицу в формуле выше) плюс дополнительно нестандартные вещественные части вида $M^r_{n-k}$, где
из $k$ прямых слагаемых $\mathbb{P}(\tau)$ можно выбрать одну пару, две пары и т.д. при условии равенства всех $c_i$; максимальное возможное число пар равно $[\frac{k}{2}]$. Отсюда следует, что число таких топологических типов дается формулой
$$
\sum_{k=0}^{n-1} (1 + [\frac{k}{2}]) = n + \sum_{k=0}^{n-1} [\frac{k}{2}] = n + [\frac{n}{2}] \cdot [\frac{n-1}{2}],
$$
что и требовалось.

При этом необходимо отметить, что полученная нижняя оценка далека от того, чтобы быть эффективной: например,  даже для простого примера $k=2$ и  частного выбора пары отображений моментов $\mu_0, \mu_1$ и их значений мы не исследовали лагранжевой геометрии раздутия $M_{n-2}$ в полном объеме.  Эта задача интересна сама по себе, поскольку
бирациональные преобразования алгебраических многообразий являются основными источниками построения новых многообразий или классификации уже известных; однако
связи лагранжевой геометрии исходного многообразия и его раздутия вдоль какого - то подмногообразия практически не исследованы (кроме достаточно простого набора примеров,
как лагранжева бутылка Клейна в поверхности Хирцебруха). Мы надеемся вернуться к этой теме в следующей работе.

$$
$$

{\bf Литература:}

[1] А. Миронов, {\it О новых примерах гамильтоново-минимальных и минимальных лагранжевых подмногообразий в $\mathbb{C}^n$ и $\mathbb{C}\mathbb{P}^n$}, Матем. сб., 195:1 (2004), 89–102;

[2] Н. Тюрин, {\it Лагранжевы циклы Миронова в алгебраических многообразиях},  Матем. сб., 212:3 (2021),  128–138;

[3] Н. Тюрин, {\it Примеры циклов Миронова в многообразиях Грассмана}, Сиб. матем. журн., 62:2 (2021),  457–465;

[4] P. Griffits, J. Harris, {\it Principles of algebraic geometry}, NY, Wiley, 1978;

[5]  M. Audin, {\it Torus action on symplectic manifolds}, 2nd rev. ed., Progr. Math.,93,Birkhäuser Verlag, Basel, 2004.

\end{document}